\documentclass[12pt]{amsart}
\usepackage{epsfig}
\newcommand{\ql}{\mathbb{C}_t[l,l^{-1},m,m^{-1}]}

\title{Skein modules and the noncommutative
torus}
\author{Charles Frohman}
\address{Department of Mathematics, University of Iowa, Iowa City, IA
52242}
\email{frohman@math.uiowa.edu}\author{R{\u{a}}zvan Gelca}
\address{Department of Mathematics, 
University of Michigan, Ann Arbor, MI 48109 and Institute of Mathematics
of the Romanian Academy, Bucharest, Romania}
\email{rgelca@math.lsa.umich.edu}

\newtheorem{theorem}{Theorem}

\newtheorem{lemma}{Lemma}
\newtheorem{corr}{Corollary}
\newtheorem{rem}{Remark}

\newcommand{\lcr}{\raisebox{-5pt}{\mbox{}\hspace{1pt}
                  \epsfig{file=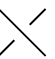}\hspace{1pt}\mbox{}}}
\newcommand{\ift}{\raisebox{-5pt}{\mbox{}\hspace{1pt}
                  \epsfig{file=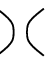}\hspace{1pt}\mbox{}}}
\newcommand{\zer}{\raisebox{-5pt}{\mbox{}\hspace{1pt}
                  \epsfig{file=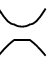}\hspace{1pt}\mbox{}}}

\begin{document}
\maketitle
\section{Introduction}

This paper introduces a new direction in the study of skein modules.
The Kauffman bracket \cite{Kau} is a knot invariant associated to quantum field
theory. The noncommutative torus is an algebra of functions that appears
in noncommutative geometry \cite{Conn}. In this paper we explicate the relationship
between the two.

When the variable of the Kauffman bracket is $-1$, the Kauffman 
bracket skein algebra of the 2-dimensional torus is isomorphic to
the algebra of $SL_2{\mathbb C}$-characters of the fundamental group of
the torus. You can think of these as  a subalgebra of the 
algebra of continuous functions on the torus. For an
arbitrary value of the variable, the Kauffman bracket
skein algebra of the torus can be viewed as a deformation of this 
particular subalgebra. Similarly, the noncommutative torus is a deformation of the
algebra of functions on the torus. The main result of the 
paper states that the Kauffman bracket skein algebra of the torus is 
isomorphic to a  subalgebra of the noncommutative torus. That is, the two
algebras arise from the same deformation.

The functions we are working with are in fact trigonometric functions,
and hence iterative techniques for dealing with Chebyshev polynomials
are a central technique for establishing the results here.
Their presence in this context is natural if one
thinks of the relation between trigonometric functions and
quantum physics.

Although a presentation of the Kauffman bracket skein algebra of the 
torus appeared before \cite{BulPrz}, the multiplicative structure of this algebra
remained mysterious.
Chebyshev polynomials with variables simple closed curves on the torus 
enable us to give a complete description of the
multiplication operation, by the {\em product-to-sum} formula
given below.

Some applications of the approach follow. First, we
analize the structure of the Kauffman bracket skein module of
the {\em solid} torus as a module over the Kauffman 
bracket skein algebra of the torus. Then, we give a short
algebraic proof of the result of Hoste and Przytycki describing
the Kauffman bracket skein module of a lens space. Finally, we show 
how to write, in terms of generators, the element of the 
Kauffman bracket skein algebra of the torus obtained by placing
a Jones-Wenzl idempotent on a simple closed curve.

\section{Skein Modules}

Throughout this paper $t$ will denote a fixed complex number.
A {\em framed link} in an orientable manifold $M$ is a 
disjoint union of annuli. In the case when the manifold can be written
as the product of a surface and an interval, framed links
will be identified with curves, using the convention that
the annulus is parallel to the surface (i.e. we consider the
blackboard framing). 
 Let $\mathcal{L}$ denote the set of equivalence
classes of framed links in $M$ modulo isotopy, including the empty link.

Consider the vector space,  $\mathbb{C} \mathcal{L}$ with  basis $\mathcal{L}$.
Define $S(M)$ to be the smallest subspace of $\mathbb{C} \mathcal{L}$
containing all expressions of the form
$\displaystyle{\lcr-t\zer-t^{-1}\ift}$
and 
$\bigcirc+t^2+t^{-2}$,
where the framed links in each expression are identical outside the
balls pictured in the diagrams. The {\em Kauffman bracket skein module} $K_t(M)$ is
the quotient
\[ \mathbb{C} \mathcal{L} / S(M). \]

Skein modules were introduced by Przytycki \cite{P} as a way to extend
the new knot polynomials of the 1980's to knots and links in arbitrary
3-manifolds. They have since become central in the theory of invariants
of 3-manifolds. The idea that they could be  used to quantize algebras
of functions on surfaces is due to Turaev \cite{T}. They were then used as
a tool for constructing quantum invariants by Lickorish \cite{Lic}, Kauffman and Lins
\cite{KL},
Blanchet, Habegger, Masbaum, and  Vogel \cite{bhmv}, Roberts \cite{Rob} and Gelca \cite{Gel}. 
Finally, the connection
between skein modules and characters of the fundamental group of the underlying
manifold was explained by Bullock \cite{B}, Przytycki and Sikora \cite{PS} and
Sikora \cite{S}. The connection between skein algebras and the algebras of observables
arising in lattice gauge field theory has been studied by Bullock, Frohman and
Kania-Bartoszy\'{n}ska \cite{BFK1}. There are also higher skein modules that
were introduced in \cite{BFK2}.

The Kauffman bracket skein module of the cylinder over a torus has 
a multiplicative structure, induced by the topological operation
of gluing one cylinder on top of the other. The product $\alpha * \beta$
is the result of laying $\alpha$ over $\beta$. This multiplication 
makes $K_t(T^2 \times I)$ into an algebra, which we will call
the skein algebra of the torus. 

The skein module of a manifold
that has a torus  boundary has a left $K_t(T^2 \times I)$-module structure 
induced by  gluing  the zero end of the  cylinder over a torus to that boundary 
component. In particular, this is true for $K_t(S^1 \times {\mathbb D}^2)$,
the Kauffman bracket skein module of the  solid torus.   
Note that $S^1 \times {\mathbb D}^2$ is homeomorphic with
the cylinder over an annulus, hence $K_t(S^1 \times {\mathbb D}^2)$
is itself an algebra. However, the algebra structure of the
skein module of the solid torus is not related to the algebra structure
of the skein module of the cylinder over the torus.
In fact $K_t(S^1 \times {\mathbb D}^2)$ is isomorphic to
${\mathbb C}[X]$ under the isomorphism that takes the
simple closed curve $\alpha $, which runs once around the torus, into the
variable $X$. Consequently, a basis of $K_t(S^1 \times {\mathbb D}^2)$
as a ${\mathbb C}$-vector space is given by the elements 
$\alpha ^n$.

\section{The Noncommutative Torus}

The {\em noncommutative torus} is a "virtual" geometric
space whose algebra of functions is a certain
deformation of the algebra of continuous functions
on the classical torus. One usually identifies the noncommutative
torus with its associated algebra of functions.

The most natural way in which the noncommutative torus 
arises is by exponentiating the Heisenberg non-commutation relation
$pq-qp=hI$. One then obtains an algebra generated 
by two unitary operators $u$ and $v$ which satisfy
$uv=\lambda vu$ where $\lambda \in{\mathbb C}$ is some constant.
The noncommutative torus is the the closure
of this algebra in a certain  $C^*$-norm.

As Rieffel \cite{Rief} pointed out, the noncommutative torus can
be obtained as a strict deformation quantization of the 
algebra of continuous functions on the torus in the following 
way. Let $t$ be the deformation parameter (denoted this way to
be consistent with the rest of the paper). For the space of
Laurent polynomials of two variables ${\mathbb C}[l,l^{-1}, m,m^{-1}]$
(here $l=exp(2\pi ix)$ and $m=exp(2\pi i y)$ are  the ``longitude'' and
the ``meridian'' of the torus), one considers the basis over ${\mathbb C}$
given by the vectors $e_{p,q}=t^{-pq}l^pm^q$. Define the multiplication
$*$, which depends on the parameter $t$, by
\[ e_{p.q}*e_{r,s}=t^{|^{p,q}_{t,s}|} e_{p+r,q+s}.\]
The space of Laurent polynomials becomes a noncommutative algebra
which we denote by $\ql$. This is the algebra of 
Laurent polynomials on the noncommutative torus. 
In order to construct the algebra 
$A_\theta$ of continuous functions on the
noncommutative torus, 
where $t=e^{2\pi i \theta}$, one considers the left regular
representation of $\ql$ on $L^2({\mathbb T}^2)$ induced by 
this product, and takes the norm closure in the operator 
norm defined by this representation. Let us mention  
that the above construction corresponds to the 
deformation of the usual product of functions in the direction of
the Poisson bracket associated to the symplectic form
$\theta dx \wedge dy$. In the physical setting mentioned at
the beginning, the unitary operators are $u=e_{1,0}$ and 
$v=e_{0,1}$. 

There is a large body of literature devoted to the algebra
$A_\theta$. In the case when $\theta$ is irrational, this 
algebra is called the irrational rotation algebra, and has appeared
in the works of operator theorists. It has been shown that
$A_\theta $ is the $C^*$-algebra naturally associated
to the Kronecker foliation of the torus $dy=\theta dx$ \cite{Conn}.
Also Weinstein explained how $A_\theta$ can be obtained through 
a geometric quantization procedure applied to the groupoid
of this foliation \cite{Wei}. 

In the present paper we are interested only in the algebra of Laurent
polynomials on the noncommutative torus. Consider the algebra morphism
\begin{eqnarray*}
{\Theta} :\ql\rightarrow\ql, \; {\Theta} (e_{p,q})=e_{-p,-q}
\end{eqnarray*}
and let  $\ql^{{\Theta}}$ be its invariant part.
Note that the algebra $\ql^{{\Theta}}$ is spanned
 by the elements $e_{p,q}+e_{-p,-q}$, $p,q\in{\mathbb Z}$.
In the next section, we will show that this algebra 
has a significant role in the study of invariants of knots. As $\Theta$
has order two, its only eigenvalues are $1$ and $-1$. The algebra $\ql$
then splits into the direct sum of its symmetric part and its antisymmetric
part with respect to $\Theta$. The subalgebra $\ql^{{\Theta}}$  is the symmetric part.

\section{The Isomorphism}

In this section we will prove that  the Kauffman bracket 
skein algebra of the torus can be embedded in the 
noncommutative torus. More precisely, we will prove that 
$K_t({\mathbb T}^2\times [0,1])$
is isomorphic to the algebra 
$\ql^{\Theta}$ defined in the previous section. 
The proof  is based on a multiplication formula,
which is the object of Theorem 1, and which is important
in its own respect. This formula  describes explicitly the multiplication
in the  Kauffman bracket  skein algebra of the torus.
Let us point out that a presentation of this algebra was 
given in \cite{BulPrz}. 
The elements that enable us to clear the picture and
obtain a neat, compact formula for the multiplication 
are {\em Chebyshev polynomials}. 

For two integers $m,n$ we denote by $\mbox{gcd}(m,n)$ their greatest
common divisor, with the convention $\mbox{gcd}(0,0)=0$.
We denote by $T_n$ the $n$-th Chebyshev polynomial,
defined recursively by $T_0=2$, $T_1=x$ and
$T_{n+1}=T_n\cdot T_1-T_{n-1}$. 

For $p,q$ relatively prime and $n\geq 0$, 
we denote by $(p,q)$ the
$(p,q)$-curve on the torus. For $(p,q)$ not necessarily
relatively prime, we define
\[
(p,q)_T=T_{\mbox{gcd}(p,q)}\left(\left(\frac{p}{\mbox{gcd}(p,q)},\frac{q}{\mbox{gcd}(p,q)}\right)
\right)
\]
which is  the  element 
 of $K_t({\mathbb T}^2\times I)$ obtained by 
replacing the variable of the Chebyshev polynomial by the 
curve on the torus.

 For $\sum \alpha_iD_i$ and $\sum \beta
_jD'_j$ two elements of the Kauffman module of the torus,
written as an algebraic combination of diagrams, we 
define their {\em intersection number} to be $max _{i,j}
D_i\cdot D'_j$, where $D_i\cdot D'_j$ is the geometric
intersection number of the diagrams $D_i$, and $D'_j$.

\begin{rem} For $m,n>0$ and $\mbox{gcd}(p,q)=1$, $\mbox{gcd}(r,s)=1$,
the geometric intersection number of $T_{n}(p,q)$ and $T_m(r,s)$ is
the absolute value of 
$mn|^{pq}_{rs}|$, where $|^{pq}_{rs}|$ is the determinant.\end{rem}

\begin{theorem} (the product-to-sum formula)
 For any  integers $p,q,r,s$ one has 
\[
(p,q)_T*(r,s)_T=t^{|^{pq}_{rs}|}
(p+q, r+s)_T+
t^{|^{pq}_{rs}|}
(p-q,r-s)_T.
\]
\end{theorem}

\proof  The proof will be by induction
on the 
intersection number of $(p,q)_T$ and $(r,s)_T$.
If the intersection number is $0$ or $\pm1$, the
relation obviously holds, by one application of the 
skein relation. 
The case $p=q=0$ or $r=s=0$ is also trivial.

\medskip

{\em Case 1.} $\mbox{gcd}(p,q)=\mbox{gcd}(r,s)=1$.

We must show that  
\[
(p,q)*(r,s)=t^{|^{pr}_{qs}|}
(p+r,q+s)_T+
t^{-|^{pr}_{qs}|}
(p-r,q-s)_T.
\]

By applying a homeomorphism of the torus, this
can be transformed in the equivalent identity:
\[
(p,q)*(0,1)=t^p
(p-1,q)_T+
t^{-p}
(p+1,q)_T
\]
with $0\leq q<p$.

If $p=1,2$, or if $q=0$, the latter equality is obvious.
To prove it for $p\geq 3$ we  use the following result.

\begin{lemma} Given $p\geq 3$, $0<q<p$ with $\mbox{gcd}(p,q)=1$ there 
exist integers $u,v,w,z$ satisfying $u+w=p$, $v+z=q$,
$|^{uv}_{wz}|=\pm 1$, $0<w<p$,
$0<u<p-1$, $0<v,z$. 
\end{lemma}

\proof  The equation
\begin{eqnarray*}
uz+vw=1
\end{eqnarray*}
can be rewritten as $u(q-v)-v(p-u)=1$ or
$uq-vp=1$. From the general theory of linear Diophantine
 equations it follows that there exists a solution 
$(u,v)$ with $0<u<p$ and
$0<v<q$. Let $w=p-u$ and 
$z=q-v$. If $u=p-1$ exchange $u$ and $w$, respectively
$v$ and $z$. \qed

Returning to the proof of the theorem, the relations $u+w=p$, $v+z=q$
and $|^{uv}_{wz}|=\pm 1$, together with the skein
relation imply
\[
(p,q)=t^{-|^{uv}_{wz}|} (u,v)*(w,z)-
t^{-2|^{uv}_{wz}|} (u-w,v-z).
\]

Since $|^{wz}_{0\, 1}|=w<p$ and $|^{u\, v}_{wz+1}|=|^{uv}_{wz}|+u=
u\pm 1<p$ we can apply the induction hypothesis to write 


\[
(p,q)*(0,1)=t^{-|^{uv}_{wz}|} (u,v)*(w,z)*(0,1)-
t^{-2|^{uv}_{wz}|} (u-w,v-z)*(0,1)=\]
\[ t^{-|^{uv}_{wz}|} (u,v)*[t^w
(w,z+1)_T+
t^{-w}
(w,z-1)_T]-\]
\[  t^{-2|^{uv}_{wz}|}[t^{u-w} 
(u-w,v-z+1)_T+
 t^{w-u}
(u-w,v-z-1)_T]=\]
 \[ t^{w+u} 
(u+w,v+z+1)_T+
 t^{-2|^{uv}_{wz}|+w-u} 
(u-w,v-z-1)_T+\]
\[ t^{-2|^{uv}_{wz}|-w+u} 
(u-w,v-z+1)_T+
 t^{-w-u} 
(u+w,v+z-1)_T-\]
\[  t^{-2|^{uv}_{wz}|+u-w} 
(u-w,v-z+1)_T+
 t^{-2|^{uv}_{wz}|-u+w} 
(u-w,v-z-1)_T=\]
\[  t^p(p,q+1)_T+
t^{-p} 
(p,q-1)_T.
\]

\medskip

{\em Case 2.} One of $\mbox{gcd}(p,q)$ or $\mbox{gcd}(r,s)$ is greater
then $1$.

Assume
that $\mbox{gcd}(p,q)\geq 2$, and let  $n=\mbox{gcd}(p,q)$, 
$p'=p/n$, $q'=q/n$. Then, an induction on $n$ gives  
\begin{eqnarray*}
& (p,q)_T*(r,s)_T=
T_{n}(p',q')*(r,s)_T=\\
& T_{n-1}(p',q')*(p',q')*(r,s)_T-T_{n-2}(p',q')*(r,s)_T=\\
 & T_{n-1}(p',q')*( t^{|^{p'r}_{q's}|}
(p'+r,q'+s)_T
+ t^{-|^{p'r}_{q's}|}
(p'-r,q'-s)_T)-\\
 & t^{(n-2)|^{p'r}_{q's}|}
((n-2)p'+r,(n-2)q'+s)_T+\\
&  t^{-(n-2)|^{p'r}_{q's}|}
((n-2)p'-r,(n-2)q'-s)_T=\\
 & t^{|^{p'r}_{q's}|+(n-1)|^{p'r}_{q's}|}
(np'+r,nq'+s)_T+
 t^{-|^{p'r}_{q's}|-(n-1)|^{p'r}_{q's}|}
(np'-r,nq'-s)_T=\\
& t^{|^{pr}_{qs}|}
(p+r,q+s)_T+
t^{-|^{pr}_{qs}|}
(p-r,q-s)_T
\end{eqnarray*}
and the theorem is proved.
\qed

\begin{theorem} There exists an isomorphism of
algebras 
\[\phi : K_t({\mathbb T}^2\times I)\rightarrow \ql^{\Theta}\]
determined by  

\[\phi ((p,q)_T)=e_{(p,q)}+e_{(-p,-q)}, \: p,q\in {\mathbb{Z}}.\]

\end{theorem}

\proof  The fact that the map is a morphism follows from
Theorem 1 and the fact that
\[
(e_{p,q}+e_{-p,-q})*(e_{r,s}+e_{-r,-s})=
t^{|^{pr}_{qs}|}
(e_{p+r,q+s}+e_{p-r,q-s})+\]
\[t^{-|^{pr}_{qs}|}
(e_{p-r,q-s}+e_{p+r,q+s}).
\]

As ${\mathbb{C}}$-vector spaces  the two algebras have the basis
$(p,q)_T$, $p\in {\mathbb{Z}}_+$, $q\in {\mathbb{Z}}$, respectively
$e_{p,q}+e_{-p,-q}$, $p\in {\mathbb{Z}}_+$, $q\in {\mathbb{Z}}$,
which proves that the map is an isomorphism. \qed

\section{The Solid Torus}

In this section we explain how to  obtain the 
skein module of the solid torus from the skein algebra of
the torus, and explicate its module structure.

 The solid torus is obtained as a quotient 
 of ${\mathbb T}^2\times [0,1]$ hence 
$K_t(S^1 \times {\mathbb D}^2)$ is obtained by factoring
$K_t({\mathbb T}^2\times [0,1])$. As mentioned before
 $K_t({\mathbb T}^2\times [0,1])$
acts on the left on the Kauffman bracket skein module of the solid torus
by the gluing map,
so the latter is a 
$K_t({\mathbb T}^2\times [0,1])$-module. Hence the the skein
module of the solid torus is the quotient of the 
skein algebra of the torus by a left ideal.
The  basis for  $K_t(S^1 \times {\mathbb D}^2)$ as a ${\mathbb C}$-vector space
  is given by $\{{\alpha }^n\}_n$, and these elements
are the images of  $(0,1)^n$, $n\geq 0$ through the
quotient map. However, for a better understanding of the module structure,
 it is better to work with the basis $\{ \alpha  _n\}_n$,  
$\alpha _n=T_n(\alpha)$. 
We denote by $\cdot $ the left action of 
$K_t({\mathbb T}^2\times [0,1])$ on the skein module of the
solid torus.

Let  ${\mathcal I}$ be the left ideal that 
is the kernel of the epimorphism
\begin{eqnarray*}
\pi :K_t({\mathbb T}^2\times [0,1])\rightarrow K_t({\mathbb D}^2\times 
{\mathbb T}).
\end{eqnarray*}
We want to show that 
$(0,1)+t^2+t^{-2}$ and $(1,1)+t^{-3}(1,0)$ form a 
minimal set of generators for ${\mathcal I}$. For this let ${\mathcal J}$
be the ideal generated by these two elements.

\begin{lemma}  Every element in $K_t({\mathbb T}^2\times [0,1])$ is of the 
form 
$P((0,1))+u $ where $P\in {\mathbb{C}}[X]$ and $u$ is in the
left ideal ${\mathcal J}$.
\end{lemma}

\proof  Since as a vector space $K_t({\mathbb T}^2\times [0,1])$ is
spanned by  $(p,q)_T$, $p,q\in {\mathbb{Z}}$, it suffices
to prove the statement for elements of this form.

If $p=0$,  since $(0,q)_T$ is a polynomial in $(0,1)$, 
\begin{eqnarray*}
(0,q)_t=a*((0,1)+t^2+t^{-2})+b, \: a\in K_t({\mathbb T}^2\times [0,1]),
 b\in {\mathbb{C}}. 
\end{eqnarray*}

If $p=1$, then from Theorem 1 we get 
\[
(1,q)_T=t^{-q+1}(1,1)*(0,q-1)_T+t^{-2q+2}(1,q-2)_T \]
and the previous argument together with an induction on $q$
show that there exist $u, u '$ in ${\mathcal J}$, 
$c\in {\mathbb{C}}$ and $P$ a
 polynomial such that
\[
(1,q)_T=u +c(1,1)+u '+P(1,0)=\]
\[u +u '+
c((1,1)+t^{-3})(1,0)-t^{-3}(1,0)+P(1,0)\]
which proves this case, too.

For $p\geq 2$ Theorem 1 gives
\[
(p,q)_T=t^{-q}(1,0)*(p-1, q)_T+t^{-2q}(p-2,q)_T 
\]
and an induction on $p$ gives the desired conclusion.

Finally, the case $p<0$ follows from
\[
(p,q)_T=t^{pq}(0,q)_T*(p,0)_T-t^{2pq}(-p,q)_T.
 \] \qed

\begin{theorem}
${\mathcal I}={\mathcal J}$. 
\end{theorem}

\proof  It is easy to see that $(0,1)+t^2+t^{-2}$ is in ${\mathcal I}$.
 On the other hand, in the solid torus, $(1,1)$ has
framing $-1$, so $(1,1)=-t^{-3}(1,0)$, from which it follows
that the second generator of ${\mathcal J}$ is in ${\mathcal I}$ as well,
hence ${\mathcal J}\subset \mathcal{I}$.

Since the restriction of  $\pi$ to the subring
of $K_t({\mathbb T}^2\times [0,1])$ generated by $(1,0)$,
previous lemma shows that ${\mathcal I}\subset {\mathcal J}$,
and the theorem is proved. \qed

Let us now describe the action of 
$K_t({\mathbb T}^2\times [0,1])$ on $K_t(S^1 \times {\mathbb D}^2)$. 
Define $x_{p,q}$ in $K_t(S^1 \times {\mathbb D}^2)$
by $x_{p,q}=t^{-pq}y_{p,q}$ where $y_{p,q}$ satisfies
\begin{eqnarray*}
& y_{p,q}=\alpha \cdot  y_{p-1,q}-y_{p-2,q}, \\
& y_{0,q}=(-t^2)^q+(-t^{-2})^q,\\
& y_{1,q}=(-t^{-3})^q(1,0).
\end{eqnarray*}

\begin{lemma}
The element  $x_{p,q}$ is the image of $(p,q)_T$ in 
$K_t(S^1 \times {\mathbb D}^2)$.
\end{lemma}

\proof 
Since
\begin{eqnarray*}
(1,0)*(p,q)_T=t^q(p+1,q)_T+t^{-q}(p-1,q)_T
\end{eqnarray*}
it follows that 
\begin{eqnarray*}
x_{p+1,q}=t^{-q}(1,0)\cdot x_{p,q}-t^{-2q}x_{p,q}
\end{eqnarray*}
from which the desired reccurence is obtained by multiplication
by $t^{(p+1)q}$. The initial condition is  a consequence of the
properties of Chebyshev polynomials.\qed

\begin{corr}
For any two integers $p$ and $q$, the element $x_{p,q} $ is a
polynomial of  $|p|$-th degree in $\alpha$. 
\end{corr}

\begin{theorem} The action of the skein algebra of the
torus on the skein module of the solid torus is given
by
\begin{eqnarray*}
(p,q)_T\cdot \alpha _n=t^{p-nq}x_{n+p,q}+t^{-p+nq}x_{p-n,q}.
\end{eqnarray*}
\end{theorem}

\proof
The theorem is a consequence of Lemma 3 and of
\begin{eqnarray*}
(p,q)_T*(n,0)_T=t^{p-nq}(n+p,q)_T+t^{-p+nq}(p-n,q)_T.
\end{eqnarray*}
\qed

\section{Lens Spaces}

In this section we will give an 
alternate short proof of the result of Hoste and Przytycki 
\cite{HosPrz} which describes the 
structure of the Kauffman bracket skein module of 
a lens space.

Let $L(p,q)$ be a lens space ($p,q\neq 0$) and
let \[ \begin{pmatrix} a & p \\  b & q \end{pmatrix}\] 
be the gluing
matrix of the two tori. Since the gluing map
reverses orientation, the determinant of this
matrix is $-1$. Any link in the lens space can be 
pushed off the cores of the two tori , thus 
\begin{eqnarray*}
K_t(L(p,q))=K_t(S^1\times {\mathbb D})\bigotimes _{K_t
({\mathbb T}^2\times [0,1])}K_t(S^1\times {\mathbb D}),
\end{eqnarray*}
where the tensor product structure is defined
by
\[x_{m,n}\otimes 1=1\otimes x_{am+pn, bm+nq}.\]
Note in particular that $K_t(L(p,q))$ is 
spanned by the elements $1\otimes 1$, $ 1\otimes \alpha$,
$1\otimes \alpha ^2$, $\cdots$. We will  prove that
in fact $K_t(L(p,q))$ is spanned by $1\otimes 1$, $1\otimes \alpha$,
$\cdots $, $1\otimes \alpha ^{[\frac{p}{2}]}$. To this end let
$V$ be the span of these elements. 

We start by noting that Theorem 3  implies that
\[ 1\otimes x_{p,q}=((0,1)\cdot 1)\otimes 1=(-t^2-t^{-2})\otimes 1
\]
and 
\[1\otimes x_{p+a, q+b}=((1,1)\cdot 1)\otimes 1=
-t^{-3}((1,0)\cdot 1)\otimes 1=-t^{-3}\otimes x_{a,b}.
\]

\begin{lemma}
For every $k\in {\mathbb Z}$, there exists a constant 
$c_k\in {\mathbb C}$ such that for all $u\in K_t(S^1\times
{\mathbb D})$ one has the identity
\begin{eqnarray*}
1\otimes ((a+kp, b+kq)\cdot u)=c_k\otimes ((a,b)\cdot u).
\end{eqnarray*}
\end{lemma}

\proof The property is true for $k=0$ and $k=1$. Since
by Theorem $1$
\begin{eqnarray*}
& 1\otimes ((a+kp, b+kq)\cdot u)=\\
& t^{aq-bp}\otimes
((p,q)*((k-1)p+a, (k-1)q+b)\cdot u)-\\ & t^{2(aq-bp)}\otimes
(((k-2)p+a, (k-2)q+b)\cdot u)=\\
& t^{aq-bp}(-t^2-t^{-2})\otimes ((k-1)p+a, (k-1)q+b)\cdot u)-\\
 & t^{2(aq-bp)}\otimes (((k-2)p+a, (k-2)q+b))\cdot u).
\end{eqnarray*}
the property follows by induction on $k$. 
\qed

\begin{lemma}
For every $m,k\in {\mathbb Z}$, one has $1\otimes x_{ma+kp,mb+kq}\in
V$.
\end{lemma}

\proof We will induct on $m$. For $m=0,1$ the
property is true, as a consequence of Lemma 4 and the fact that
$1\otimes x_{kp,kq}=(-t^2-t^{-2})^k\otimes 1$. 

Let $k_0$ be the integer that minimizes  the absolute value of
$ma+k_0p$. Clearly this minimum is at most $[\frac{p}{2}]$, hence
$x_{ma+k_0p,mb+k_0q}$ is in $V$. On the other hand, by using Theorem 1, 
we get that, for an arbitrary $k$,
\begin{eqnarray*}
& x_{ma+kp, mb+kq} =t^{-(m-1)k-mk_0}\otimes ((a+(k-k_0)p, b+(k-k_0)q)*\\
& ((m-1)a+k_0p, (m-1)b+k_0q)\cdot 1)-\\ &t^{-2(m-1)k-2mk_0}\otimes
x_{(m-2)a-(k-2k_0)p, (m-2)b-(k-2k_0)q}.
\end{eqnarray*}
>From Lemma 4 and Theorem 1 it follows that
\begin{eqnarray*}
& 1\otimes (a+(k-k_0)p, b+(k-k_0)q)*\\
& ((m-1)a+k_0p, (m-1)b+k_0q)\cdot 1)=\\
& c_k\otimes ((a,b)*((m-1)a+k_0p, (m-1)b+k_0q)=\\
& c_kt^{-1}x_{ma+k_0p, mq+k_0q}+c_ktx_{(m-2)q+k_0p,(m-2)b+k_0q}.
\end{eqnarray*}
 So, from the induction hypothesis and the fact that
$x_{ma+k_0p,mb+k_0q}\in V$ we get that 
$x_{ma+kp, mb+kq}\in V$, which completes the induction.\qed

\begin{theorem} (Theorem 4. in \cite{HosPrz}) The space
$K_t(L(p,q))$ is spanned by $1\otimes 1$, $1\otimes \alpha$,
$\cdots $, $1\otimes \alpha ^{[\frac{p}{2}]}$.
\end{theorem}

\proof Every natural number $n$ can be written in the form
$ma+kp$. From  Corollary 1 it follows that 
\[
1\otimes x_{ma+kp, mb+kq}=c\otimes \alpha ^{ma+kp}+1\otimes f(\alpha)
\]
where $f$ is a polynomial of degree strictly less than $ma+kp$,
and $c$ is a nonzero constant, hence  by applying Lemma 5 and inducing on
$n$ we deduce that $\alpha ^n\in V$ for every $n$, which proves the
theorem.\qed 

\section{Jones-Wenzl idempotents}

Jones-Wenzl idempotents appeared for the first time in the 
study of operator algebras, but they are best known to 
topologists because of their use in the construction of
topological quantum field theories (\cite{Lic, bhmv, Rob, Gel}).
By placing Jones-Wenzl idempotents on simple closed curves on
a torus one obtains certain elements of the skein algebra of
the torus. We show below how one can make use of the
embedding of this algebra in the noncommutative torus
to give a pleasing formula for these skeins.

Fix $r$ an integer greater than $1$, and let $t=e^{i\pi/2r}$. For
a positive integer $n$, the $n$-th Jones-Wenzl idempotent $f^{(n)}$ lives
in the Temperley-Lieb algebra $TL_n$, which, let us remind, is the
algebra of diagrams of non-intersecting strands joining 
$2n$ points on the boundary of a rectangle, with multiplication
defined by juxtaposition of rectangles.
Jones-Wenzl idempotents are denoted by empty boxes, and are defined
inductively as in Fig. 1. Here the convention is that a number $k$
written next to a strand means $k$ parallel strands, and 
$\Delta _k=(-1)^k[k+1]$, where 
$[k+1]$ is the quantized integer 
$\frac{t^{2k+2}-t^{-2k-2}}{t^2-t^{-2}}$.
Recall that Jones-Wenzl idempotents are defined only for
$n=0,1,\cdots , r-2$. 

\begin{figure}[htbp]
\centering
\leavevmode
\epsfxsize=4in
\epsfysize=.8in
\epsfbox{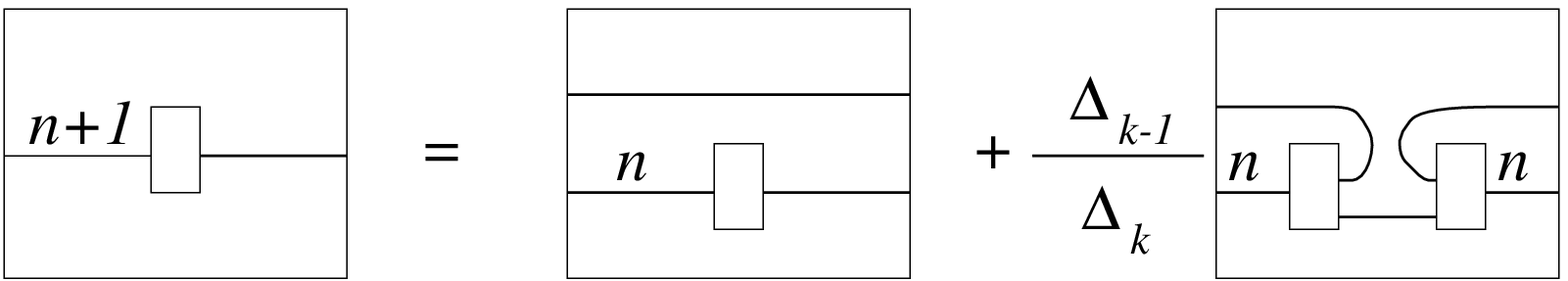}

Figure 1.  
\end{figure}

We will denote by $(p,q)_{JW}$ the element of $K_t({\mathbb T}^2\times I)$
obtained by taking gcd$(p,q)$ parallel copies of the 
$(\frac{p}{\mbox{\small 
gcd}(p,q)},(\frac{q}{\mbox{\small gcd}(p,q)})$-curve and 
inserting on them the gcd$(p,q)$-th Jones-Wenzl idempotent.

\begin{theorem}
If $p$ and $q$ are relatively prime and $n$ is a positive integer less 
than $r-2$, then
\begin{eqnarray*}
& (np, nq)_{JW}=(np,nq)_T+((n-2)p, (n-2)q)_T+\cdots \\
& +((n-2[\frac{n}{2}])p,
(n-2[\frac{n}{2}])q)_T
\end{eqnarray*}
\end{theorem}

\begin{figure}[htbp]
\centering
\leavevmode
\epsfxsize=3in
\epsfysize=1in
\epsfbox{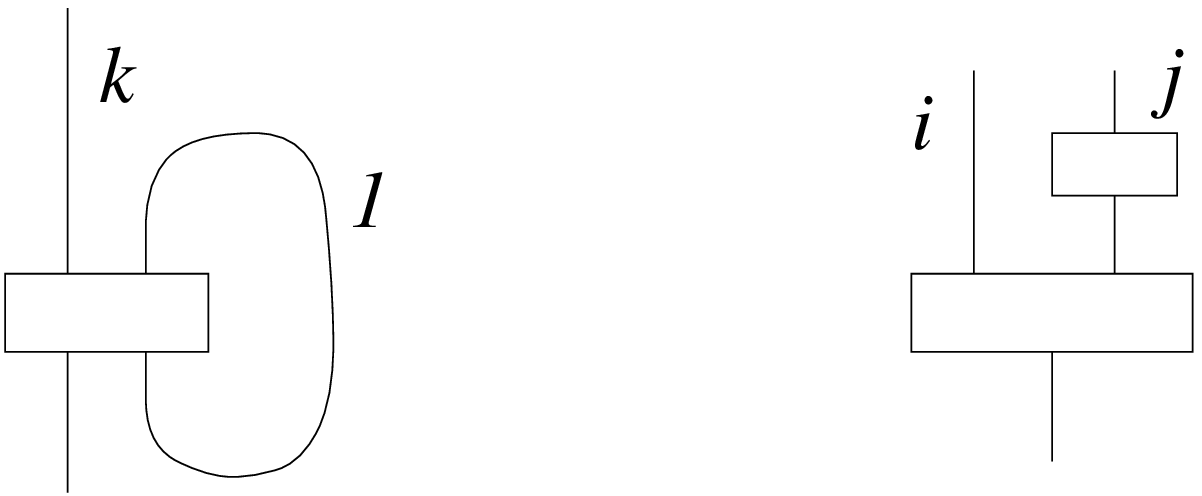}

Figure 2. 
\end{figure}

\proof By using the well known identities from Fig. 2
we deduce that the following recurrence relation holds
\begin{eqnarray*}
(np,nq)_{JW}=(p,q)*((n-1)p,(n-1)q)_{JW}-((n-2)p, (n-2)q)_{JW}.
\end{eqnarray*}
By Theorem 2,  the image of the $(p,q)$-curve in the noncommutative 
torus is $e_{p,q}+e_{-p,-q}$. In the noncommutative torus
the 
elements \[e_{mp,mq}+e_{(m-2)p,(m-2)q}+\cdots +e_{-(m-2)p,-(m-2)q}+e_{-mp,-mq}\]
with $m\geq 0$, satisfy the same recurrence relation as $(np,nq)_{JW}$.
Since the image of $f^{(0)}$ is $1$, and the image of
$(p,q)_{JW}$ is $e_{p,q}+e_{-p,-q}$ ($p$ and $q$ are relatively
prime),  the conclusion
of the theorem follows by induction.\qed

Note that both $(np,nq)_{T}$ and $(np,nq)_{JW}$ satisfy the recurrence relation
of Chebyshev polynomials. Also the term corresponding to
$n=1$ is the same in both cases. However, the term corresponding
to $n=0$ is $2$ in the first case, and $1$ in the second, which
makes them so much different.

It is a well known fact that the trace of a Jones-Wenzl idempotent is
\begin{eqnarray*}
& (-1)^n [n+1]=(-1)^n\frac{t^{2n+2}-t^{-2n-2}}{t^2-t^{-2}}=\\
& (-1)^n
(t^{2n}+t^{2n-2}+\cdots +t^{-2n+2}+t^{-2n})
\end{eqnarray*}
 As the theorem and its proof show the elements
$(np,nq)_{JW}$ are obtained by replacing 
$(-t^2)^k+(-t^{-2})^k$ in this formula with $(kp,kq)_T$.

\end{document}